# A Proof of Collatz Conjecture Based on a New Tree Topology

**Hassan Rezai Soleymanpour**

E-mail address: rezai.s.hassan@gmail.com

**Abstract** – Consider a finite positive integer. If it is even, divide it by 2, and if it is odd, multiply it by 3 and add 1. This will give you a new integer. Following the procedure for the new integer, you will receive another integer. Repeat the steps, and after a few repetitions, you will finally reach 1. Collatz conjecture states that the final integer in the mentioned process will always be 1, no matter what integer it starts with. Although the procedure of the conjecture is easy to describe, its correctness has not yet been confirmed. This article proves the conjecture by introducing a tree topology of it. Given the proposed tree, we can prove that all integers are uniquely distributed on the tree, and there is no cycle other than 1-2-1. We also see how every trajectory ends in 1.

## I- The Collatz function

Let $f : \mathbb{N} \to \mathbb{N}$ be the function of Collatz conjecture that is defined as follows [1,2]:

$$f(d) = \begin{cases} \dfrac{3d+1}{2} & \text{if } d \text{ is odd} \\ \dfrac{d}{2} & \text{if } d \text{ is even} \end{cases} \qquad (1)$$

The Collatz process starts with a finite integer and is repeated using the return value of the function as the input argument for the next iteration. The conjecture states that for every integer $d \in \mathbb{N}$ there is always an integer $t \in \mathbb{N}$ so that $f_t(d) = 1$ where $f_t$ represents the $t^{\text{th}}$ iteration of the function $f$ .

## II- An overview of the proposed tree topology

To get acquainted with some definitions, an overview of the proposed tree is first provided in this section. The proposed tree topology of Collatz conjecture is shown in Fig. 1. Because it expands to infinity, the tree cannot be displayed for all integers, so (...) is placed to indicate that the tree can be continued from that point. The idea behind this topology is to arrange all integers in such a way that for every integer, all possible integers that reach it in the Collatz process are placed to the right or top of it. In this way, if an integer is repeated anywhere in the tree, all the integers attached to it must be repeated at the same point.

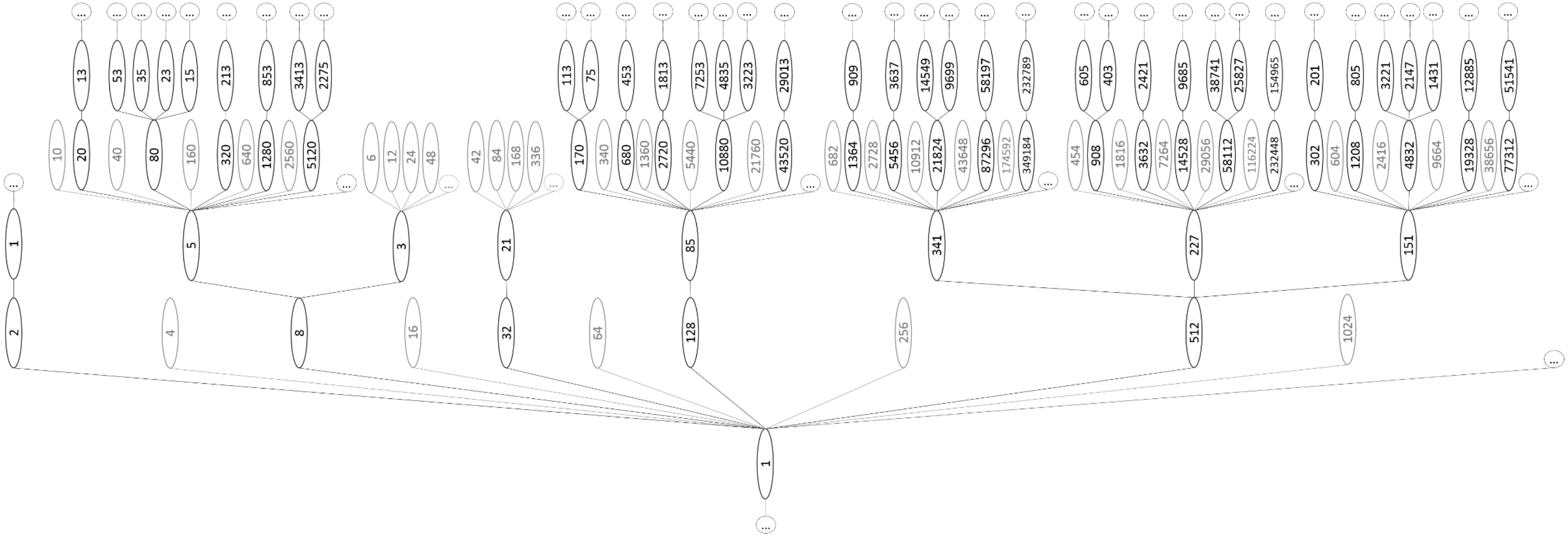

*Figure 1-The proposed tree of Collatz conjecture*

The proposed tree consists of an infinite number of sub-trees. Every sub-tree is defined by a parent (the node) directly below some children (the branches). The children of a sub-tree are all the integers that reach the parent during the Collatz process. A sub-tree is called odd or even, depending on what the parent is. An odd sub-tree consists of an odd parent with an infinite number of even children, while an even sub-tree consists of an even parent with a certain number of odd children. All sub-trees are hierarchically connected, so that the children of a sub-tree are the parents for their overhead sub-trees while the parent of the sub-tree itself is a child of its underneath sub-tree. Fig. 2 shows an example of two sub-trees. The left sub-tree (sub-tree (a)) is an even sub-tree because its parent, 80, is even, and the children, 53, 35, 23, and 15 are odd. On the other hand, the right sub-tree (sub-tree (b)) is an odd sub-tree because the parent, 85, is odd, and the children, "170, 340, 680, 1360, …" are all even. During the Collatz process, all the children of a sub-tree finally reach their parent. For example, 15, 23, 35, and 53 reach 80, and "170, 340, 680, 1360, …" reach 85. Note that 15, 23, 35, and 53 are the only odd integers that reach 80, and "170, 340, 680, 1360, …" are the only even integers that reach 85.

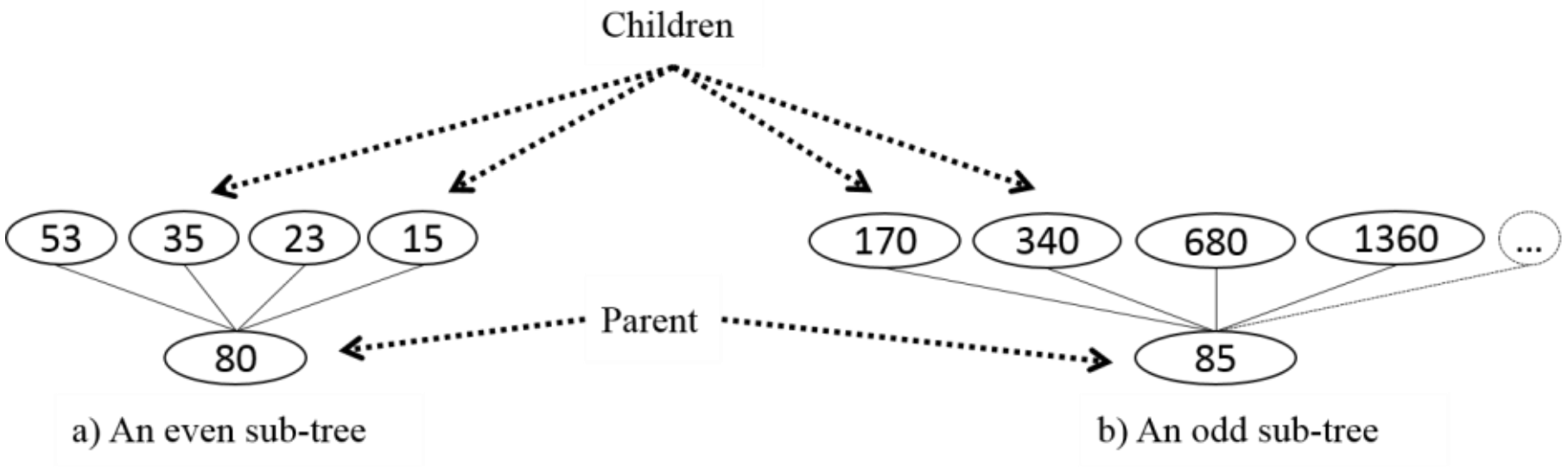

*Figure 2-An Example of even and odd sub-trees*

## III- Even sub-trees

Let $D \in \mathbb{N}$ be an odd integer by which the return values of the Collatz function are odd for the next $n-1$ iterations, and the return value of the $n^{\text{th}}$ iteration is even. So, the even integer $f_n(D)$ after iteration $n$ can be calculated as follows:

$$\begin{aligned} f_1(D) &= \frac{3}{2}D + \frac{1}{2} \\ f_2(D) &= \frac{3}{2}\left(\frac{3}{2}D + \frac{1}{2}\right) + \frac{1}{2} = \left(\frac{3}{2}\right)^2 D + \frac{1}{2}\left(\frac{3}{2} + 1\right) \\ &\vdots \\ f_n(D) &= \left(\frac{3}{2}\right)^n D + \frac{1}{2}\sum_{i=1}^{n}\left(\frac{3}{2}\right)^{i-1} \end{aligned} \tag{2}$$

Substituting $\sum_{i=1}^{n}\left(\frac{3}{2}\right)^{i-1} = \dfrac{\left(\frac{3}{2}\right)^n - 1}{\frac{3}{2} - 1}$, Eq. (2) is simplified as follows:

$$f_n(D) = \left(\frac{3}{2}\right)^n (D+1) - 1 \tag{3}$$

The odd integer $D$ can be written as follows:

$$D(n,k) = 2^n k - 1 \qquad n \in \mathbb{N}, k \in \text{odd} \tag{4}$$

where $k$ is a positive odd integer. In this way Eq. (3) is simplified as follows:

$$f_n(D) = T(n,k) = 3^n k - 1 \tag{5}$$

It means that the odd integer $D(n,k) = 2^n k - 1$ reaches the even integer $T(k,n) = 3^n k - 1$ after $n$ consecutive iterations. For example, 31 can be written as $31 = D(5,1) = 2^5 \times 1 - 1$ with $k = 1$ and $n = 5$. Therefore, it reaches $T(5,1) = 3^5 \times 1 - 1 = 242$ after 5 iterations.

Let $k$ be represented by its prime factors according to the fundamental theorem of arithmetic as follows:

$$k = 3^{n_1} 5^{n_2} 7^{n_3} \ldots \qquad n_1, n_2, n_3, \ldots, \in \mathbb{Z}_+ \tag{6}$$

Therefore Eq. (5) is modified as follows:

$$T(n,k) = 3^n k - 1 = 3^n (3^{n_1} 5^{n_2} 7^{n_3} \ldots) - 1 \tag{7}$$

Let $h = 5^{n_2} 7^{n_3} 11^{n_4} \ldots$ be an odd integer so that $3 \nmid h$ (3 does not divide $h$ ) and $N = n + n_1$, then Eq. (7) becomes:

$$T(N,h) = 3^N h - 1 \qquad N \in \mathbb{N} \tag{8}$$

Eq. (4) is also modified accordingly as follows:

$$D(n,N,h) = 2^n 3^{(N-n)} h - 1 \qquad n,N \in \mathbb{N}, \quad n \le N \tag{9}$$

The two above equations state that there are exactly $N$ odd integers of form Eq. (9) for $1 \le n \le N$ that eventually reach the same even integer of form Eq. (8) during the Collatz process. Therefore, they can shape an even sub-tree with the $N$ odd integers as the children and the even integer as the parent. For example, let $N = 4$ and $h = 5$, so there are 4 children as follows:

$$D_1(n=1, N=4, h=5) = 2^1 \times 3^{(4-1)} \times 5 - 1 = 269$$
$$D_2(n=2, N=4, h=5) = 2^2 \times 3^{(4-2)} \times 5 - 1 = 179$$
$$D_3(n=3, N=4, h=5) = 2^3 \times 3^{(4-3)} \times 5 - 1 = 119$$
$$D_4(n=4, N=4, h=5) = 2^4 \times 3^{(4-4)} \times 5 - 1 = 79$$

That all of them eventually reach $T(N=4, h=5) = 3^4 \times 5 - 1 = 404$ . Fig. 3 shows this even sub-tree. Note that if we start the Collatz process from 79 (the child with the least value), we reach the other children one by one (i.e. $79 \to 119 \to 179 \to 269$ ), and finally the parent 404.

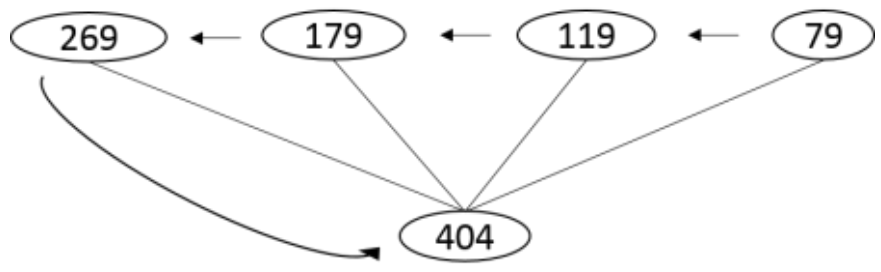

*Figure 3 –The even sub-tree when* $N = 4$ *and* $h = 5$ .

The configuration of an even sub-tree is shown in Fig. 4. The black circle represents the even parent while the white circles show the $N$ odd children. The arrows show how each member reaches the other in the Collatz process. We can show that every child reaches its left child during the Collatz process as follows: Let $D_i(i, N, h) = 2^i 3^{(N-i)} h - 1$ be a child of an even sub-tree. The left adjacent child is $D_{i-1}(i-1, N, h)$ which can be written as follows:

$$D_{i-1} = 2^{i-1} 3^{(N-(i-1))} h - 1 = \frac{3}{2}(2^i 3^{(N-i)} h - 1) + \frac{1}{2} = \frac{3}{2} D_i + \frac{1}{2} \tag{10}$$

Which means that $D_{i-1}$ is the Collatz value of $D_i$ . It is also obvious that $T > D_1 > D_2 > D_3 > ... > D_N$ .

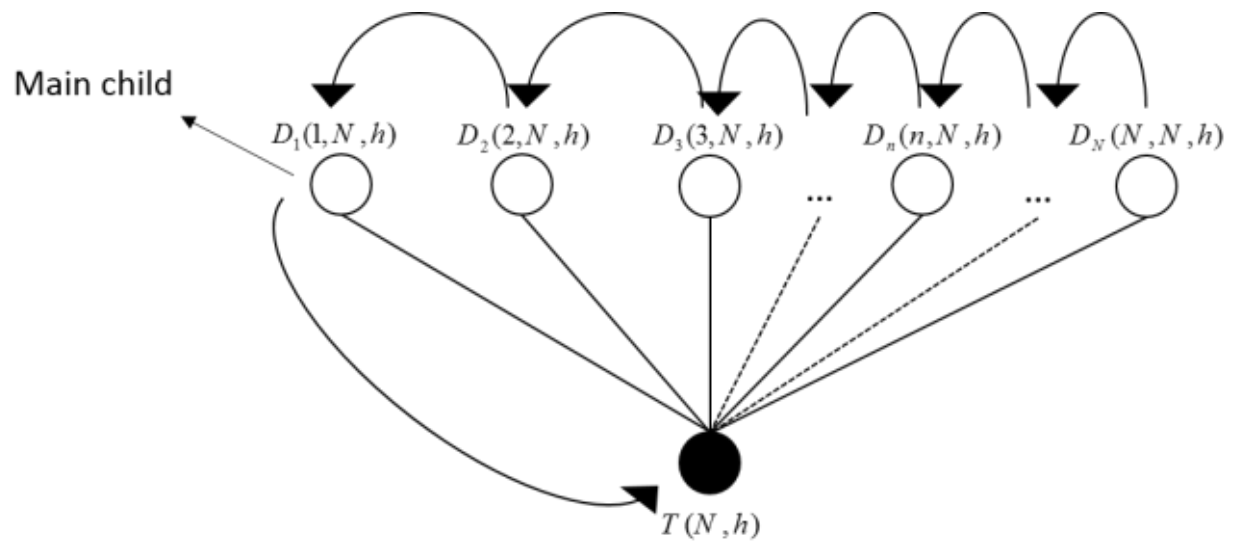

*Figure 4- the configuration of an even sub-tree*

Because it always exists in the sub-tree and is the only odd integer that reaches the parent directly, the first child of every even sub-tree (i.e. $D_1(1, N, h) = 2 \times 3^{(N-1)} h - 1$ ) is called the **main child** of the sub-tree. It means that the other children first reach it before they reach the parent. The main child also has the highest value among the other children. We will later use this concept to show how every trajectory ends in 1.

***Every even sub-tree is unique to each of its members***

According to the fundamental theorem of arithmetic, any odd child $D(n, N, h) = 2^n 3^{(N-n)} h - 1$ is unique to the set $\{n, N, h\}$, and any even parent $T(N, h) = 3^N h - 1$ is unique to the set $\{N, h\}$. In other words, by having a set of $n, N$, and $h$ we have a unique odd child, and by having a set of $N$ and $h$ we have a unique even parent. Because $1 \le n \le N$ , an even sub-tree is identified when we are given the set $\{N, h\}$. We can easily calculate the set $\{N, h\}$ from each of the members. For example, suppose that we are given 359 as a child of an even sub-tree, and we are asked to identify the other members. We have $359 = 2^3 \times 3^2 \times 5 - 1$ which means that $n = 3$ , $N = 5$ and $h = 5$ . Therefore, the sub-tree has 5 children which are $D(1,5,5) = 809$ , $D(2,5,5) = 539$ , $D(3,5,5) = 359$ , $D(4,5,5) = 239$ , and $D(5,5,5) = 159$ , and the parent is $T(5,5) = 1214$ . There is no even sub-tree other than this, in which we see any of these members. Therefore, whenever we see any of these members in the tree, we must see the other members around it, and whenever one of these members is repeated somewhere on the tree, the other members must be repeated too.

**IV- Odd sub-trees**

An odd sub-tree consists of an infinite number of even children that reach an odd parent. The children are all even integers that can be calculated by multiplying $2^b\ ,b\in\mathbb{N}$ by the parent. Unlike even sub-trees that have a unique form, there are three types of odd sub-trees according to their parents that are investigated in this section.

We observed that every odd child eventually reaches an even parent in the form of $T=3^N h-1$ . The next step in the Collatz process is to divide the even integer by 2 enough to reach another odd integer. Because $3\nmid\left(3^N h-1\right)$, if we divide $T=3^N h-1$ enough by 2, it eventually reach an odd integer of the following set:

$$\{1,5,7,11,13,17,19,...\}=\{1,7,13,...\}\bigcup\{5,11,17,...\}=\{6a-5|\,a\in\mathbb{N}\}\bigcup\{6a-1|\,a\in\mathbb{N}\}\tag{11}$$

The above set includes all odd integers except those that can be divided by 3, i.e. $\{6a-3\,|a\in\mathbb{N}\}=\{3,9,15,21,...\}$. It means that when an even integer in the form of $T=3^N h-1$ is divided enough by 2 to reach an odd integer, the odd integer does not belong to the set $(6a-3)$ . Therefore, any odd integer in the form of $(6a-3)$ has no even integer of form $T=3^N h-1$ above itself, and since all odd integers eventually reach an even integer in the form of $T=3^N h-1$ , there is no odd integer above any odd integer of form $(6a-3)$. For example, there is no odd integer above 21 in the Collatz tree, and there is no even integer of form $T=3^N h-1$ above it because $21=(6\times4-3)$ is divisible by 3. But still, we can consider every odd integer of form $(6a-3)$ as the parent of infinite even children of form $(6a-3)2^m\quad a\in\mathbb{N},m\in\mathbb{N}$ . For example, 21 is the parent of an odd sub-tree that has infinite children, i.e. "42, 84, 168, …". As a result, every odd sub-tree which has a parent of form $(6a-3)$ terminates its trajectory to higher rows in the Collatz tree. Therefore, we can call these sub-trees the *flowers* of the Collatz tree, because they are at the top of the branches.

Since the even integers of form $3^N h-1$ eventually reach odd parents of the set (11), they can be written as either $(6a-1)2^b,\,b\in\mathbb{N}$ or $(6a-5)2^b,\,b\in\mathbb{N}$. But this is not true for all values of $b$. Because $3^N h-1=3(3^{N-1}h)-1=3k-1=6i-4\quad i\in\mathbb{N}$ where $k=3^{N-1}h$ is an odd number, we can check for which values of $b$ , we have $6i-4=(6a-1)2^b$ or $6i-4=(6a-5)2^b$ .To do so, we should see for which values of $b$ we have $i=\dfrac{(6a-1)2^b+4}{6}$ or $i=\dfrac{(6a-5)2^b+4}{6}$ to be integers. For $i=\dfrac{(6a-1)2^b+4}{6}$ we have:

$$i=\frac{(6a-1)2^b+4}{6}=a2^b-\frac{2^b-4}{6}\tag{12}$$

Let us decompose $b$ to odd and even values. For even values of $b$ , we have:

$$b = 2m \ , m \in \mathbb{N} \Rightarrow i = a2^{2m} - \frac{2^{2m}-4}{6} = a4^{m} - \frac{1}{2}(\frac{4^{m}-1}{3} - 1) \tag{13}$$

$3 \mid \left(4^{m}-1\right)$ because $4^{m}-1 = (4-1)(4^{m-1}+4^{m-2}+...+1)$. Therefore, Eq. (12) is an integer for even values of $b$. For odd values of $b$, we have:

$$b = 2m-1 \ , m \in \mathbb{N} \Rightarrow i = a2^{2m-1} - \frac{2^{2m-1}-4}{6} = \frac{1}{2}a4^{m} - \frac{1}{4}(\frac{4^{m}-1}{3} - \frac{7}{3}) \tag{14}$$

Since $3 \mid \left(4^{m}-1\right)$, Eq. (12) is not an integer for odd values of $b$.

Similarly, we can do the above process for $i = \frac{(6a-5)2^{b}+4}{6}$ as follows:

$$i = \frac{(6a-5)2^{b}+4}{6} = \frac{(6a-6+1)2^{b}+4}{6} = (a-1)2^{b} + \frac{1}{2}\frac{2^{b}+4}{3} \tag{15}$$

For even values of $b$, we have:

$$b = 2m \ , m \in \mathbb{N} \Rightarrow i = (a-1)2^{2m} + \frac{1}{2}(\frac{2^{2m}+4}{3}) = (a-1)4^{m} + \frac{1}{2}(\frac{4^{m}-1}{3} + \frac{5}{3}) \tag{16}$$

Since $3 \mid \left(4^{m}-1\right)$, Eq. (15) is not an integer for even values of $b$. For odd values of $b$, we have:

$$b = 2m-1, m \in \mathbb{N} \Rightarrow i = (a-1)2^{2m-1} + \frac{1}{2}(\frac{2^{2m-1}+4}{3}) = \frac{1}{2}(a-1)4^{m} + \frac{1}{4}(\frac{4^{m}-1}{3} + 3) \tag{17}$$

Because $3 \mid \left(4^{m}-1\right)$, Eq. (15) is an integer for odd values of $b$. To summarize the results, every even integer of form $T(N,h) = 3^{N}h - 1$ can be written as either one of the following forms:

$$T(N,h) = 3^{N}h - 1 = \begin{cases} (6a-1)2^{2m} \\ or \\ (6a-5)2^{(2m-1)} \end{cases} ,a \in \mathbb{N}, m \in \mathbb{N} \tag{18}$$

For example, $1700 = T(5,7) = 3^{5} \times 7 - 1 = (6 \times 71 - 1) \times 2^{2\times1}$ which has the form of $(6a-1)2^{2m}$ where $a = 71$ and $m = 1$. As another example, $842 = T(1,281) = 3^{1} \times 281 - 1 = (6 \times 71 - 5) \times 2^{(2\times1-1)}$ which has the form of $(6a-5)2^{(2m-1)}$ where $a = 71$ and $m = 1$.

By looking at Eq. (18), we observe that there are infinite even children of form $(6a-5)2^{(2m-1)}$ for $m = 1,2,3,...,\infty$ that reach the unique odd parent $(6a-5)$, and there are infinite even children of form $(6a-1)2^{2m}$ for $m = 1,2,3,...,\infty$ that reach the unique odd parent $(6a-1)$. The even integers of forms $(6a-1)2^{(2m-1)}$ and $(6a-5)2^{2m}$ (i.e. those even integers that cannot be written as $3^{N}h - 1$ ) can also be considered as the children of the sub-trees although they do not have any odd integers directly over themselves.

*We can summarize the results as follows:*

- There are three forms of odd sub-tree according to the three parent forms: $(6a-1)$ , $(6a-3)$ and $(6a-5)$ .
- The three forms of odd parent contain all odd integers, and the children of them contain all even integers because they can be obtained by multiplying the parents by all powers of 2.
- Every parent of even sub-trees is a child of an odd sub-tree of the form either $(6a-1)$ or $(6a-5)$ and it can be written as either $(6a-1)2^{2m}$ or $(6a-5)2^{(2m-1)}$ .
- The even integers of forms $(6a-3)2^{m}$ , $(6a-1)2^{(2m-1)}$ , or $(6a-5)2^{2m}$ cannot be written as $T(N,h)=3^{N}h-1$. Therefore, they cannot be parents for any even sub-trees.
- There is no odd integer above the odd integers of form $(6a-3)$ . Therefore, every odd sub-tree of form $(6a-3)$ terminates its trajectory to higher rows. Because of this, we can call these sub-trees the flowers of the Collatz tree.

Figures 5 to 7 show the three types of odd sub-tree. The black circles indicate that the even children can be written as $T(N,h)=3^{N}h-1$ and they are the parents of their superior sub-trees. The gray circles are those even children which cannot be written as $T(N,h)=3^{N}h-1$, so they do not have upper sub-trees. The white circles represent the odd parents of the sub-trees.

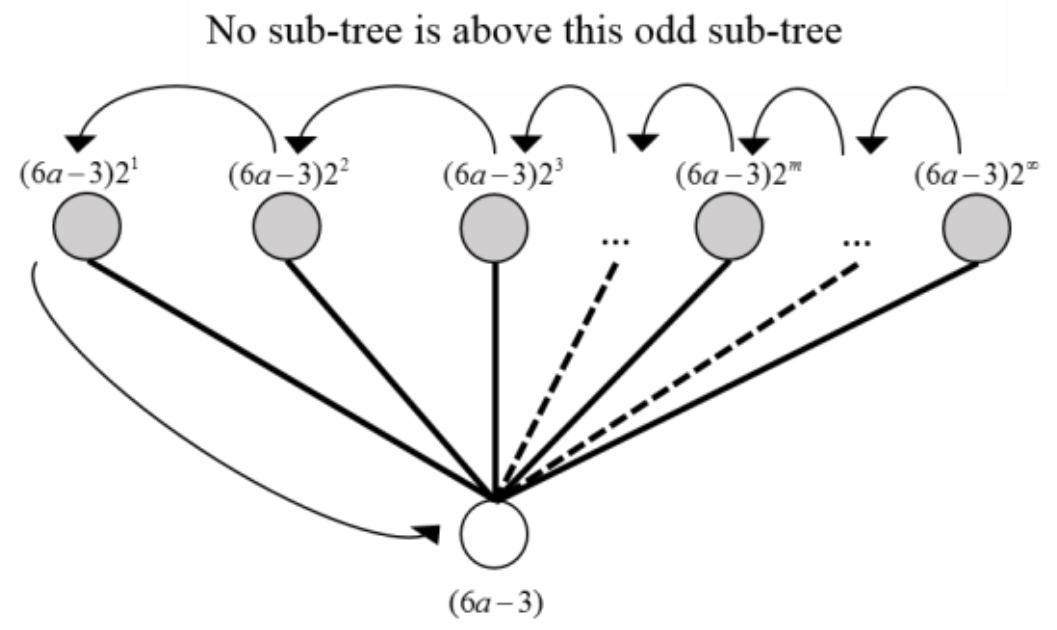

*Figure 5- the configuration of an odd sub-tree of form (6a-3)*

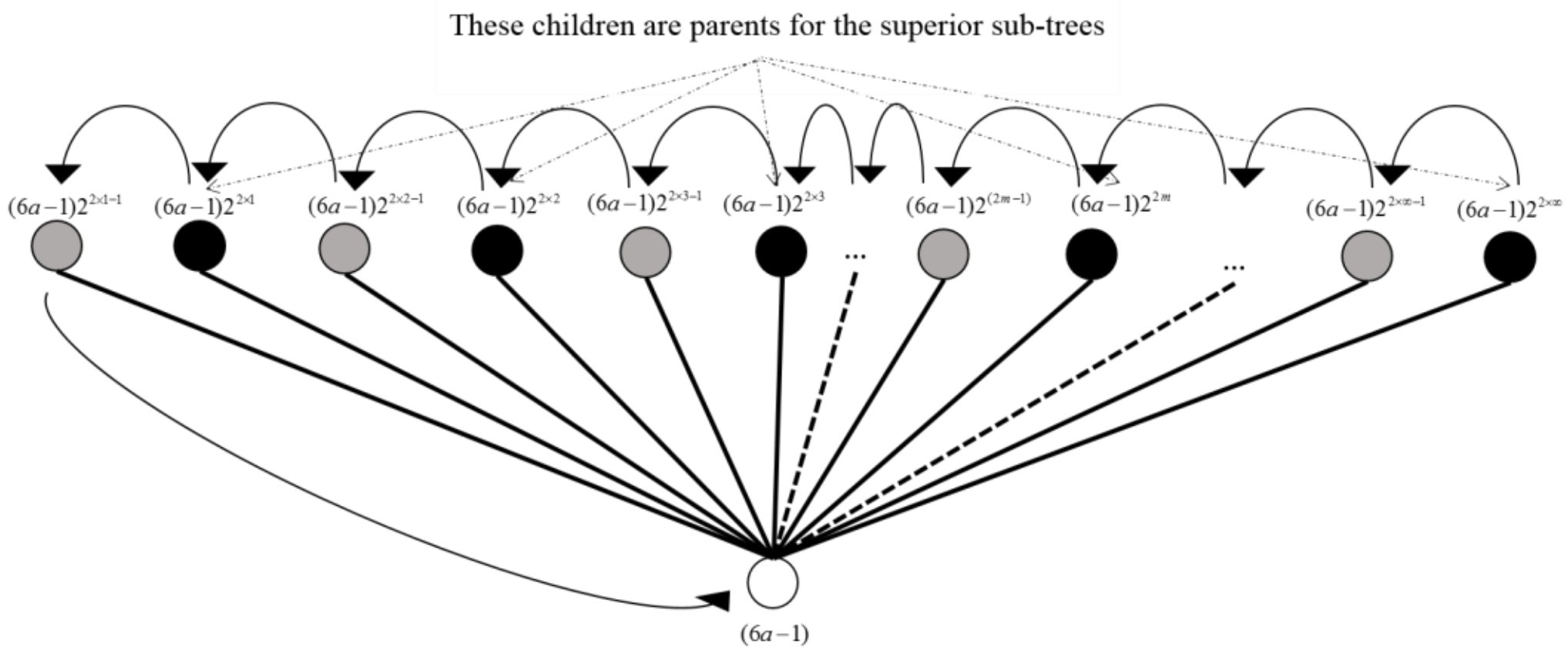

*Figure 6- the configuration of an odd sub-tree of form (6a-1)*

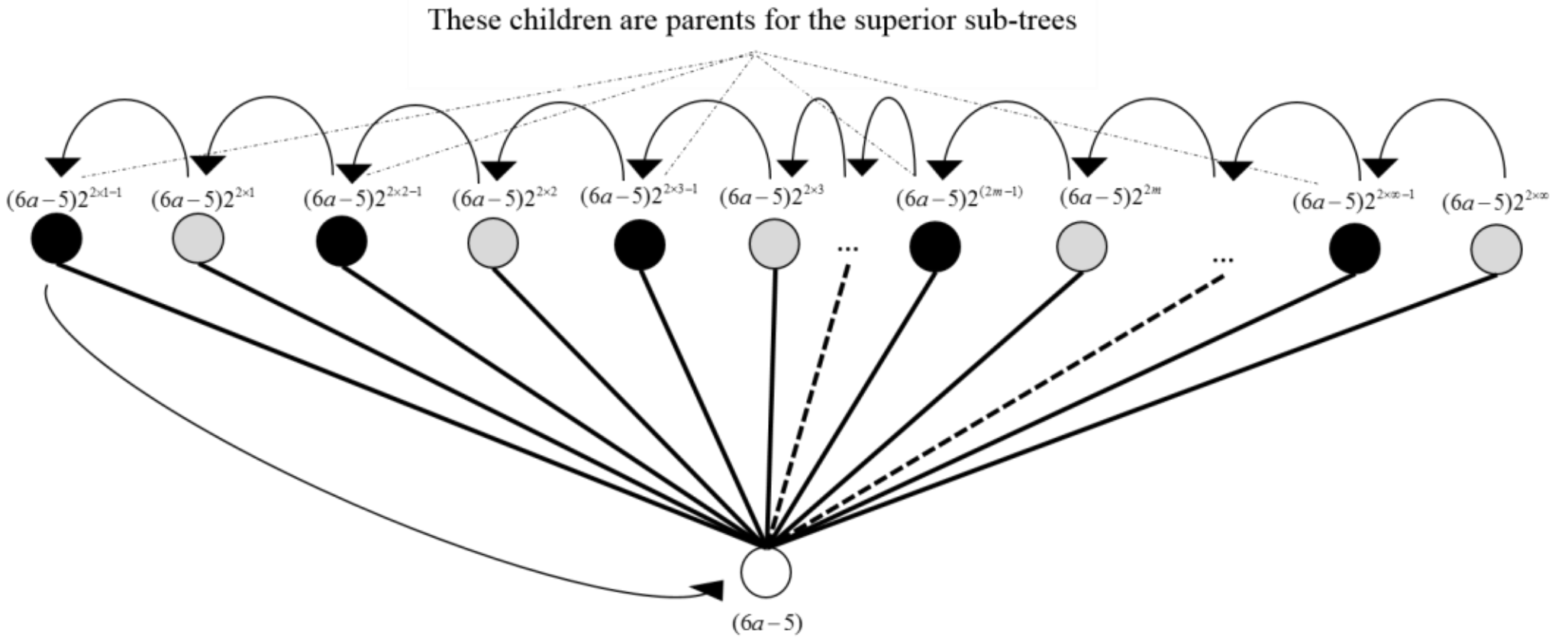

*Figure 7- the configuration of an odd sub-tree of form (6a-5)*

The arrows on the figures show the direction of movement during the Collatz process. Each child reaches its adjacent left child because it is equal to the left child multiplied by 2. The direction of movement is from right to left and then to the beneath row, similar to even sub-trees.

***Every odd sub-tree is unique to each of its members***

The parent of an odd sub-tree is an odd integer that has one of the following three forms: $(6a-1)$, $(6a-3)$, or $(6a-5)$. The children are all the even integers that can be obtained by multiplying the parent by all powers of 2. Therefore, when we are given an even child, we can obtain the parent uniquely by factorization according to the fundamental theorem of arithmetic. For example, $1152 = 2^7 \times 9$. Therefore, 1152 is the 7th child of the sub-tree and the parent is 9 which is of form $(6a-3)$ where $a = 2$. It means that 1152 is only the child of this sub-tree and, it cannot be a child of any other sub-tree. Therefore, an odd sub-tree can be identified whenever a member of it is known, and we cannot find any two different odd sub-trees that have at least a member in common.

## V- Drawing the tree

To draw the tree, we start from 1 as the parent of the first sub-tree and then continue the tree by considering the children of each sub-tree to be the parents of the overhead sub-trees. We know that 1 is an odd integer, so it is the parent of an odd sub-tree that can be written as $1 = (6\times1-5)$. Therefore, we have an odd sub-tree of type $(6a-5)$ where $a = 1$. The children of the sub-tree are $(6\times1-5)\times2^m, m \in \mathbb{N} = 2,4,8,...$, but those children that extend the tree are $(6\times1-5)2^{(2m-1)}, m \in \mathbb{N} = 2,8,32,...$. Let us continue the tree with 2, which is the first child of 1. Because 2 is an even parent, it can be written as $2 = T(1,1) = 3^1\times1-1$ which means that it has one child, and the child is $D(1,1,1) = 2^1\times3^{1-1}\times1-1 = 1$. We see that the child and parent of 2 are equal, therefore 1 and 2 are repeated infinitely. But remember that 1 has infinite children and whenever it is repeated, all of its children are repeated too, so the tree has no starting point. Therefore, the real topology of Collatz conjecture is something like Fig. 8 where the tree is repeated infinitely. But the goal is to see if all integers eventually reach 1, so we ignore the cycle of 1 and 2.

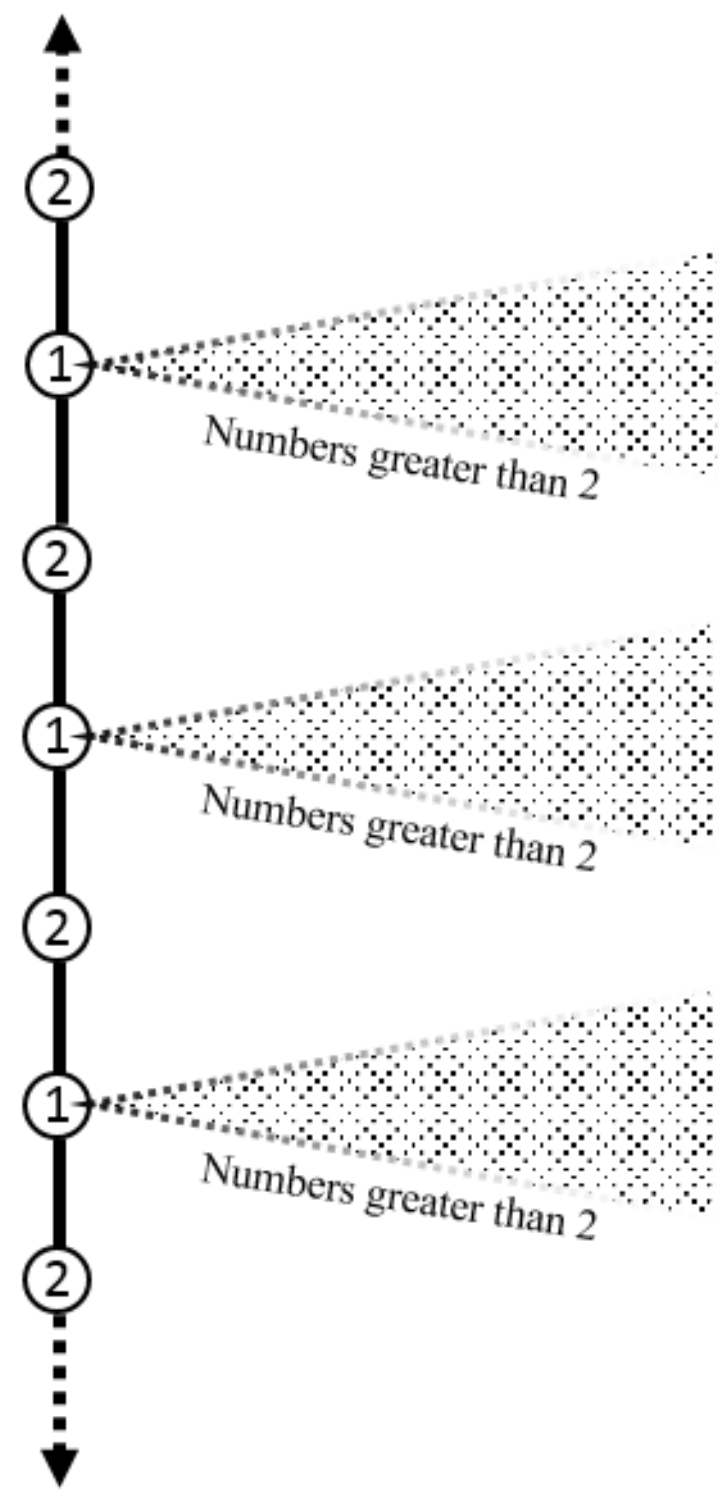

*Figure 8- The real topology of Collatz conjecture*

Let us continue with the other children of 1. The second child is 4 which is not an even parent. The third child is 8 which is an even parent so that it can be written as $8 = T(2,1) = 3^2 \times 1 - 1$. Therefore, it has 2 children which are $D(1,2,1) = 2^1 \times 3^{2-1} \times 1 - 1 = 5$ and $D(2,2,1) = 2^2 \times 3^{2-2} \times 1 - 1 = 3$. We see that 5 is an odd parent of type $(6a-1)$ with $a=1$ and it has infinite even children i.e. "10, 20, 40, …" that half of them, i.e. $(6\times1-1)2^{(2m)} = 20, 80, 320, ...$ extend the tree because they are even parents. On the other hand, 3 is an odd parent of type $(6a-3)$ with $a=1$ which has infinite even children i.e. "6, 12, 24, …" but none of them is an even parent, so there is no other sub-tree above 3. We can continue the tree from the rest of the children of 1 (i.e. "32, 128, 512, …") as well as the children of 5.

## VI- The proof of Collatz Conjecture

In the previous sections, we observed that every even integer is a child of an odd sub-tree, and every odd integer is a child of an even sub-tree. It means that every integer reaches a parent, and the parent reaches another parent. Therefore, all integers are connected, and there is no discontinuity on the tree, although we may have two or more separate trees. Because of this, when we start from an integer we reach another integer and this process is repeated forever unless we trap in a cycle. Since the number of children is always greater than or equal to 1 in even sub-trees and it is countless in odd sub-trees, the tree is growing from a root, except for when there is a cycle in the tree's body that makes the tree repeat from a point. Besides, we observed that every sub-tree is unique to each of its members in such a way that if a member is repeated somewhere in the tree, all the members attached to it, and hence all the sub-trees that are connected must be repeated. It means that we do not have a cycle of finite numbers, and if there is a cycle, a tree of numbers must be repeated. For example, we saw that when 1 and 2 are repeated, all numbers connected to them are repeated too (Fig. 8). Therefore, if there is a cycle at the tree's root, the whole tree is repeated (like Fig. 8). On the other hand, if there is a cycle in the tree's body, then the tree must be symmetrical from that point. So, if we prove that every integer has a unique place in the tree, then we prove that there is no cycle in the tree's body, and if we prove that there is only one tree with a unique root, we prove the conjecture. Let us first see how many roots we can find for the Collatz topology as follows.

***- How many separate trees exist?***

Let us suppose that there is no cycle in the tree's body. We know that the number of members in each row is bigger than its bottom row, so the tree members converge to a root. It means that when we start from a child and continue the Collatz process, it will be terminated to a root either by one of the following conditions: a) there has to be a root integer that does not have a parent below itself, and b) there has to be a cycle that traps the trajectory in a repetitive root. The first condition is not true because, for every integer, there is always a unique parent. But, as we know, the second condition is true. For example, 1 and 2 are both child and parent of each other making a cycle. We should check that if there is any root other than one so that we have separate trees. To do so, we must look for a child that reaches a parent so that the parent reaches the child again. It means that the child and its grandparent are the same. So, we must look for a 1-trivial cycle other than 1-2-1. Although it had been proved in [4] that there is no 1-trivial cycle other than 1-2-1, another proof is presented in this paper as follows.

Let $D(n, N, h)$ be an odd child of an even sub-tree that reaches the even parent $T(N, h)$. If we divide $T(N, h)$ enough by 2, we reach the parent of $T(N, h)$ which in this case is equal to the child of it i.e. $D(n, N, h)$. Therefore, we have:

$$D(n, N, h) = 2^n 3^{(N-n)} h - 1 \tag{19}$$

$$T(N, h) = 2^b D(n, N, h) = 3^N h - 1 \quad b \in \mathbb{N} \tag{20}$$

Substituting Eq. (19) in Eq. (20) and doing some rearrangements we reach the following equation:

$$(2^b 2^n - 3^n)k = 2^b - 1 \tag{21}$$

where $k = 3^{(N-n)} h$ is an odd number, and it is obvious that $2^b > \left(\frac{3}{2}\right)^n$. Let us subtract Eq. (19) from Eq. (20) as follows:

$$2^b D - D = 3^N h - 2^n 3^{(N-n)} h \Rightarrow (2^b - 1)D = k(3^n - 2^n) \tag{22}$$

In the above equation $D = 2^n k - 1 = k$ if $n = k = 1$, and they are relatively prime integers otherwise. Therefore, if $n = 1$ and $k = 1$ we have $b = 1$ and $D = 1$. Let us assume $k > 1$ or $n > 1$ we reach the following system of equations (note that when $n > 1$, $(3^n - 2^n)$ and $(2^b - 1)$ are relatively prime integers and so are $k$ and $D$ when $k > 1$):

$$\begin{cases} 2^b - 1 = k \\ D = 2^n k - 1 = 3^n - 2^n \end{cases} \tag{23}$$

Because $2^b - 1 = k$ , Eq. (21) becomes:

$$(2^b 2^n - 3^n)k = k \Rightarrow 2^{b+n} - 1 = 3^n \tag{24}$$

$b+n$ must be even because $3 \nmid (2^{b+n} - 1)$ for odd values of $b+n$ . On the other hand, If $b+n$ is even, then $n$ must be odd because $3^n$ cannot be perfect square. So, $b$ and $n$ are both odd integers. Because $b$ is odd $2^b - 1 = k$ in Eq. (23) is not divisible by 3. Therefore, $k$ is of forms $(6a-1)$ and $(6a-5)$ . If $k = 2^b - 1$ is of form $(6a-1)$ , we have: $6a - 1 = 2^b - 1 \Rightarrow 3a = 2^{b-1}$ which is not possible because of 3. Therefore, $k$ is of form $(6a-5)$ . So we have: $2^b - 1 = (6a-5) \Rightarrow 2^b = 2(3a-2)$ . Therefore, $b = a = k = 1$, and because $2^b > \left(\frac{3}{2}\right)^n$ we have $n = D = 1$.

This shows that the only odd child that can be the root of the Collatz tree is $D(1,1,1) = 1$ , and the only 1-trivial cycle is 1-2-1. So we have one tree with a unique root.

***- Is there any cycle in somewhere in the tree's body?***

We proved that the only child that can be equal to its grandparent is 1. Let $G$ be an odd integer of type either $(6a-1)$ or $(6a-5)$ . The type $(6a-3)$ is not in our interest, because there is no odd integer above it. We know that there are countless even sub-trees above $G$ which the parents of them are of type either $(6a-1)2^{2m}$ or $(6a-5)2^{(2m-1)}$ depending on $G$ . Also, every even sub-tree above $G$ has at least one child. Therefore, $G$ has an infinite number of odd grandchildren which we can consider in a block like the one shown in Fig. 9.

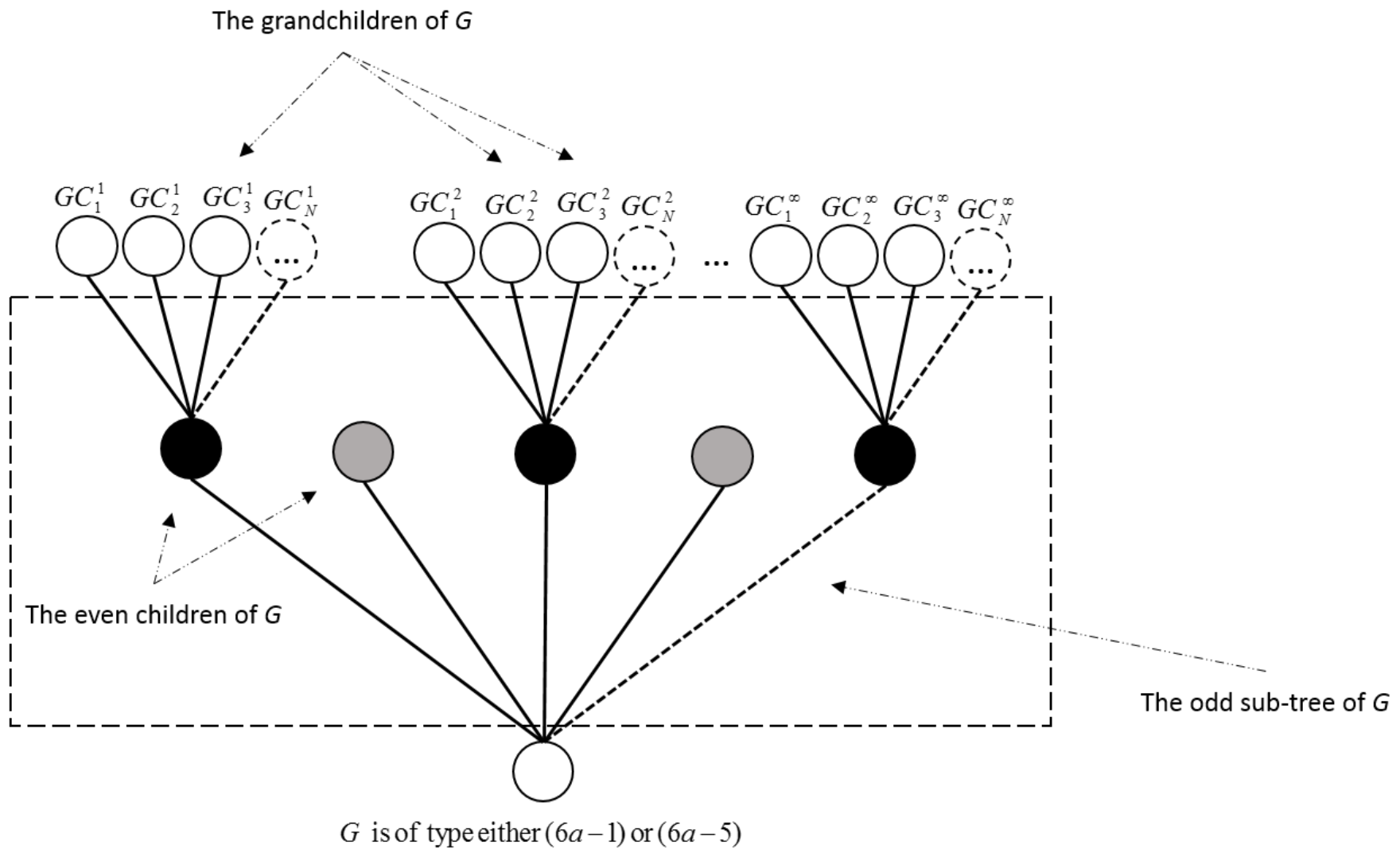

*Figure 9- The block of G to its grandchildren*

$GC_i^j$ represents a grandchild of $G$ which is the $i^{\text{th}}$ child of the $j^{\text{th}}$ productive child of $G$ . The block of $G$ to its grandchildren has the following properties:

1- There is an infinite number of grandchildren for $G$ which all of them are odd integers. Because we have all the possible even children of $G$ in the block, the grandchildren are all the possible ones for $G$ .

2- All of the grandchildren are unique to $G$ and each other because of the uniqueness of their sub-trees. It means that any of these grandchildren only belongs to this block when they have the rule of a grandchild and we have: $GC_i^j \neq GC_l^z \ \forall i, j, z, l \in \mathbb{N}$.

3- If $G > 1$, we have $GC_i^j \neq G \ \forall i, j \in \mathbb{N}$ (By the proof of the previous section).

4- If $G$ does not belong to any cycle, then all of the grandchildren do not belong to any cycle because every grandchild finally reaches the grandparent and then continues the Collatz process on the grandparent's trajectory.

We observed that the only odd integer that can be the root of the Collatz tree is 1. Also, we know that 1 has infinite grandchildren that are "1, 5, 3, 21, 85, …" that all of them except 1 do not belong to any cycle because if they are, they never reach 1 as their grandparent (Look at Fig. 10). According to the mentioned properties, the grandchildren of 5 are all unique integers that only available in the block of 5. Also, Because 5 is not in a cycle, any of its grandchildren is not in a cycle too. This concept happens for the grandchildren of "85, 341, 227, 151, …". Because "3, 21, …" are of form $(6a-3)$ they have no grandchildren. Therefore, the grandchildren of every row are all unique odd integers that are not in the following rows. Because their grandparents do not belong to any cycle, they also do not belong to any cycle. Therefore, we see that we have new unique odd integers in every row without trapping in any cycle. It means that all integers are distributed uniquely on the tree, and they finally reach the tree's root that is 1. So, the conjecture is **verified**.

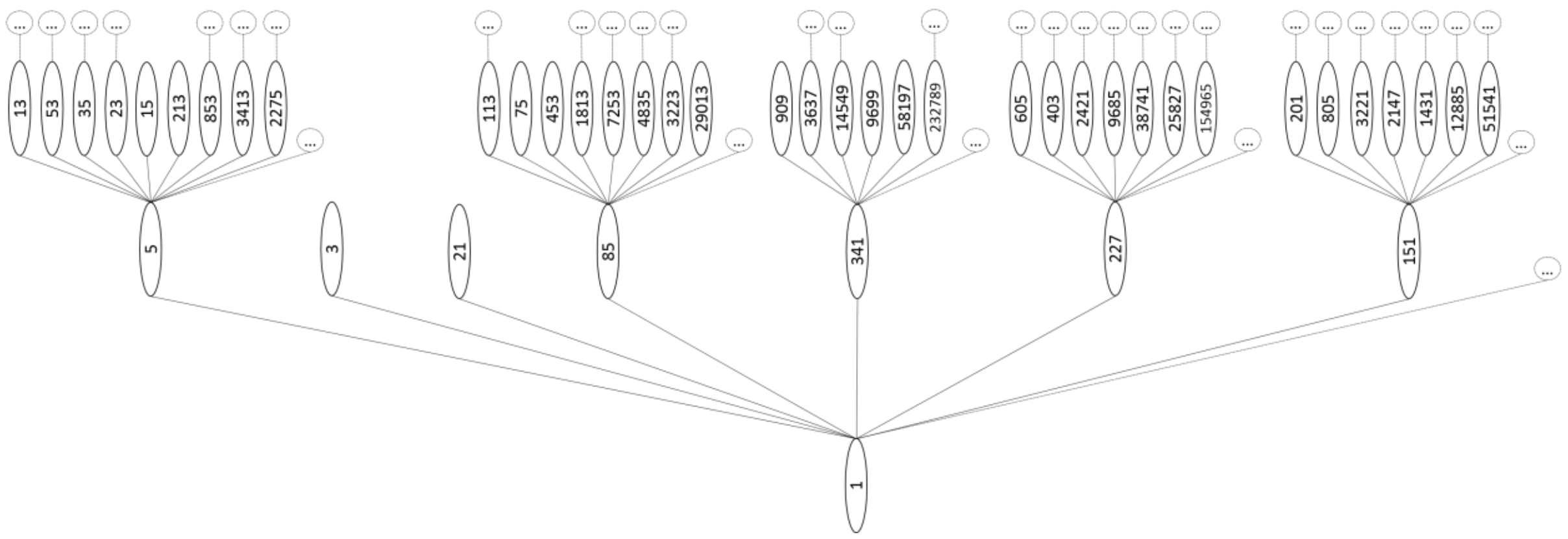

*Figure 10-The Collatz tree showing only odd integers*

## VII- How every trajectory eventually ends in 1

In section III, the concept of the main child of an even sub-tree is defined, and we observed that all odd children in the even sub-tree eventually reach the main child before they reach the parent. Let $D(1, N, h) = 2 \times 3^{(N-1)} h - 1$ be the main child of an even sub-tree. The parent of the sub-tree is $T(N, h) = 3^N h - 1$, which is itself a child of the parent of the following odd sub-tree, which we call it the grandparent (see Fig. 11). We can easily prove that the main child is always greater than its grandparent except when the main child is equal to 1.

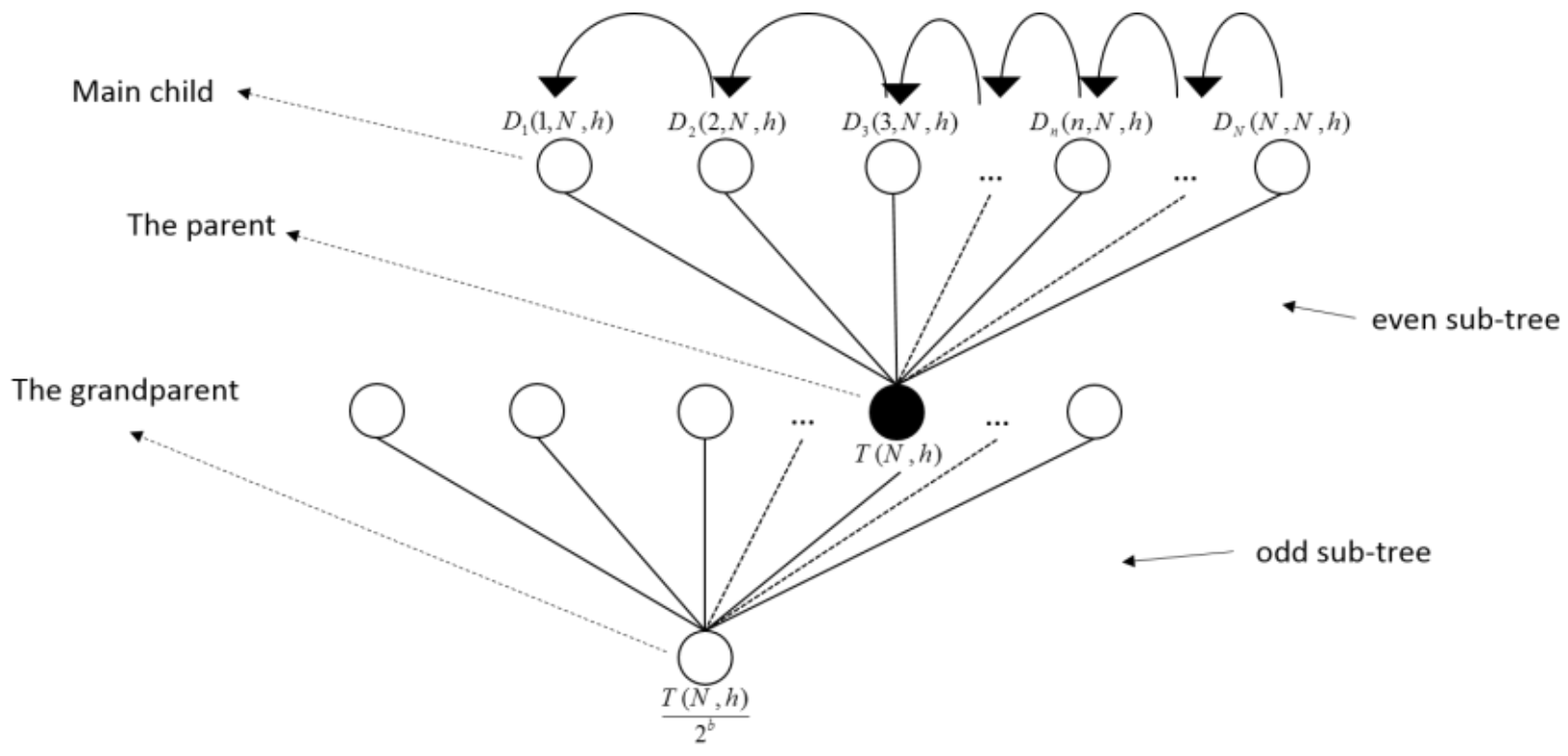

*Figure 11- A main child and its grandparent.*

We know that the grandparent is obtained by dividing $T(N,h)$ by $2^b$, $b \in \mathbb{N}$. Therefore, we claim that:

$$D(1,N,h) \geq \frac{T(N,h)}{2^b} \tag{25}$$

If the above inequality is valid for $b=1$, that is the minimum value of $b$, then it is certainly valid for the other values of $b$. So we continue with $b=1$ as follows:

$$2\times 3^{(N-1)}h-1 \geq \frac{3^N h-1}{2} \tag{26}$$

By some rearrangements, we reach the following inequality:

$$3^{N-1}h \geq 1 \tag{27}$$

This inequality states that the main child is always greater than its grandparent when $N$ or $h$ is greater than 1, and it is equal to it when $N$ and $h$ are both equal to 1, i.e. when the main child is equal to 1. For example, see Fig. 1. We see that 13, 53, 213, 853, and 3413 which are the main children of their sub-trees are all greater than their grandparent (i.e. 5), and 5 is itself greater than 1. Let $G>1$ be an odd integer of form either $(6a-1)$ or $(6a-5)$. The main child of every even sub-tree above $G$ can be calculated as follows:

$$Mc_1(G,b_1) = \frac{2^{b_1}G-1}{3} \tag{28}$$

where $Mc_1$ represents the main children of the first row above $G$ and $b_1 = 2m_1+1$, $m_1 \in \mathbb{N}$ if $G$ is of type $(6a-1)$ or $b_1 = 2m_1$, $m_1 \in \mathbb{N}$ if $G$ is of type $(6a-5)$. Let us calculate the main child of every even sub-tree of the $J^{\text{th}}$ row above $G$ as follows:

$$Mc_2(G,b_1,b_2)=\frac{2^{b_2}(\frac{2^{b_1}G-1}{3})-1}{3}=\frac{2^{b_1+b_2}G-2^{b2}-3^1}{3^2}$$

$$Mc_3(G,b_1,b_2,b_3)=\frac{2^{b_3}(\frac{2^{b_1+b_2}G-2^{b2}-3^1}{3^2})-1}{3}=\frac{2^{b_1+b_2+b_3}G-2^{b_2+b_3}-3^1 2^{b3}-3^2}{3^3} \tag{29}$$

$\vdots$

$$Mc_J(G,b_1,b_2,...,b_J)=\frac{2^{\sum_{i=1}^{J}b_i}G-\sum_{u=1}^{J}3^{u-1}2^{\sum_{i=u+1}^{J}b_i}}{3^J}$$

where $b_i=2m_i+1, m_i\in\mathbb{N}$ if $Mc_{i-1}$ is of type $(6a-1)$ or $b_i=2m_i\ , m_i\in\mathbb{N}$ if $Mc_{i-1}$ is of type $(6a-5)$. We see that $Mc_i$ is unique to $G$ and $b_1,b_2,...,b_i$ according to the fundamental theorem of arithmetic. Although the values $b_1,b_2,...,b_i$ could be selected the same for countless main children, the main child $Mc_i$ is unique to the odd integer $G$ and it eventually reaches it. It means that:

$$\forall G_1,G_2\in\{\{6a-5\}\cup\{6a-1\}\}\quad Mc_i(G_1,b_1,b_2,...,b_i)=Mc_i(G_2,b_1,b_2,...,b_i)\Rightarrow G_1=G_2. \tag{30}$$

For example, Let $G=5$, the main children of the first row above 5 can be calculated as follows: $Mc_1(5,m_1)=\frac{2^{2m_1+1}5-1}{3},\ m_1\in\mathbb{N}=13,53,213,853,...$. We know that 13 is of form $(6a-5)$. Therefore, The main children of 13 which belong to the set of main children of the second row above 5 are $Mc_2(5,m_1=1,m_2)=\frac{2^{2(1)+1+2m_2}5-2^{2m_2}-3}{3^2},\ m_2\in\mathbb{N}=17,69,277,1109,...$. Let us see the main children of the third row above 5 when $m_1=1$, and $m_2=3$ as follows: $Mc_3(5,m_1=1,m_2=3,m_3)$

$=\frac{2^{2(1)+1+2(3)+2m_3}5-2^{2(3)+2m_3}-3^1 2^{2m_3}-3^2}{3^3},\ m_2\in\mathbb{N}=369,1477,5909,23637,...$ . Fig. 12 shows these main children.

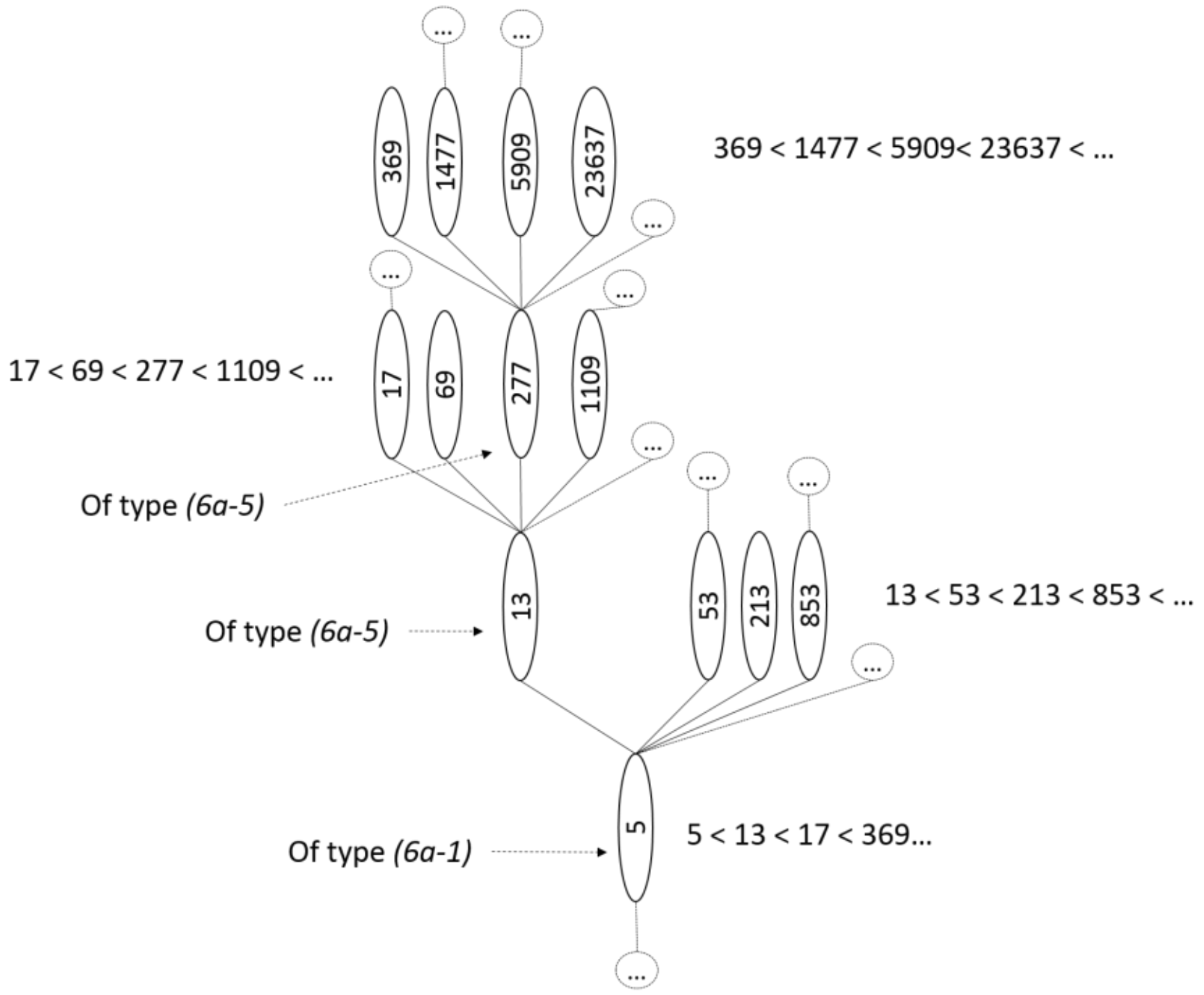

*Figure 12- the trajectories of 5*

We observe that $369 < 1477 < 5909 < 23637 < ...$ and $17 < 369$, and also $17 < 69 < 277 < 1109 < ...$ and $5 < 13 < 17$ which means that the main children are sorted from bigger to smaller when we see them from right to left, or from up to down. Also, we see that all main children are unique. We can call these branches the trajectories to 5 which are always from greater to smaller. We observe that 5, itself is a main child reaching its grandparent, 1. Also, for every main child there is a unique set of children that reach it before they reach the parent. So, the process is as follows. Let us start the Collatz process from an odd integer greater than one. If the integer is a main child it is on a trajectory to an odd integer which is not a main child (except if the odd integer is 1). If the integer is not a main child, it reaches the main child of its sub-tree. For example, 3077 is a main child and it is on the trajectory of 23 which is placed on the fifth row above 23. On the other hand, 23 is not a main child, but it reach to 53 which is the main child of its sub-tree, and this process is repeated for 53 as follows. 53 is the main child and it is on the trajectory to 1. Therefore, 3077 finally reaches 1 through the mentioned trajectories.

We can call the trajectories to 1 as the *major trajectories* of the Collatz tree because all other trajectories finally end in them. Fig. 13 shows the major trajectories. The other trajectories are started from the children that are attached to the main children of the major trajectories. This makes the tree expanding using new unique integers. Because the main children of the superior sub-trees are greater than those of the following sub-trees, we see bigger numbers when we move to higher rows of the tree. This concept is also true when we move to the right side of the tree.

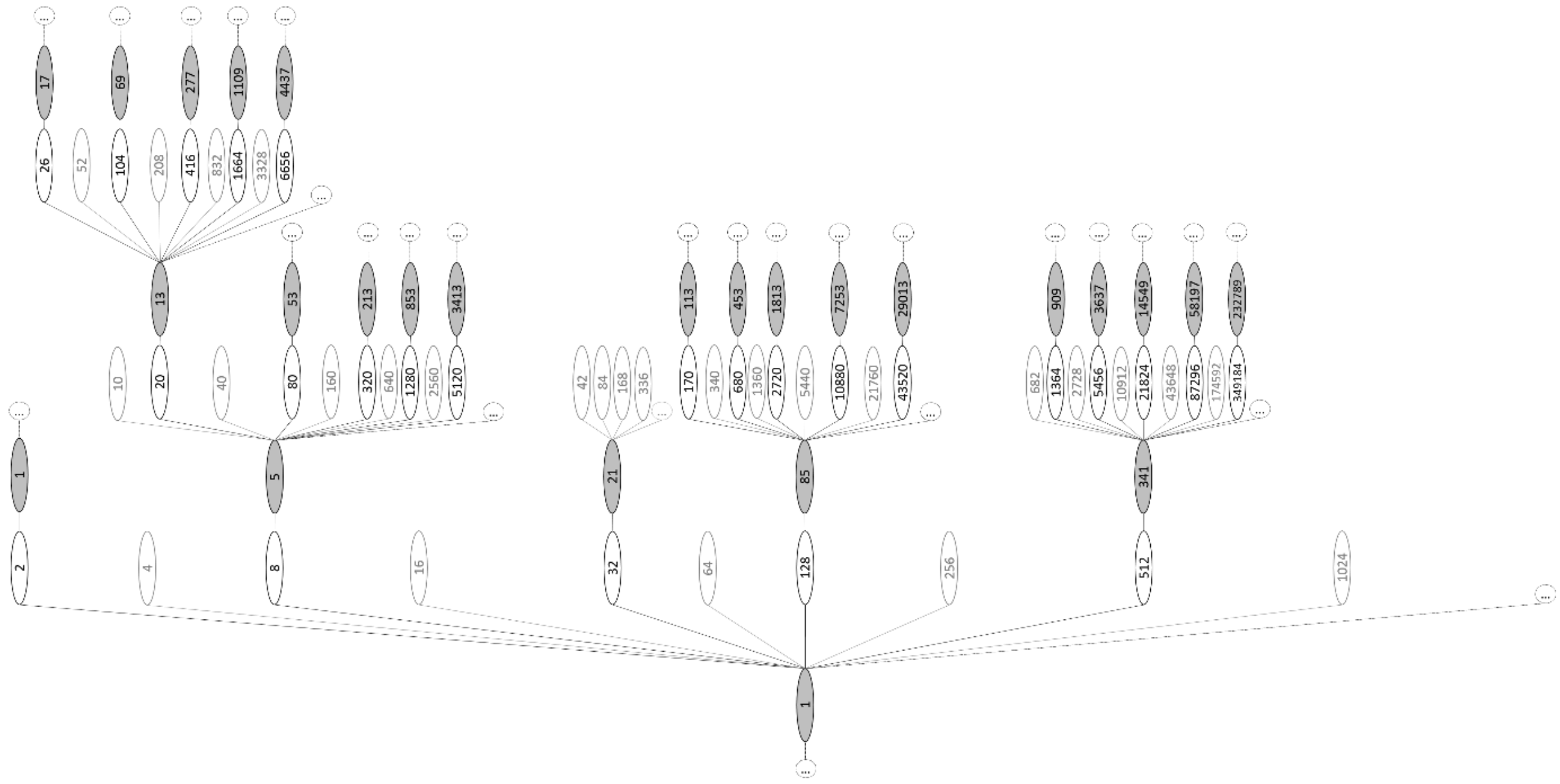

*Figure 13- The tree with only major trajectories*

## VIII- Conclusion

A new topology of Collatz conjecture is presented in this paper. We observed that the Collatz topology is a tree on which every integer has a unique place. All integers that are divisible by 3, i.e. “3, 6, 9, 12, …” shape the flowers of the tree, while the rest of them form the branches. The tree consists of an infinite number of sub-trees that are hierarchically connected. Every sub-tree is unique to its members so that whenever a member is repeated somewhere all other members must be repeated too. According to this, we prove that there is no cycle other than 1-2-1 and all integers finally reach 1, and then trap in that cycle. We also introduce the concept of “main child”, and we see that the main children are always in unique trajectories to smaller numbers. This property shows how every trajectory finally ends in 1.

## IX- References

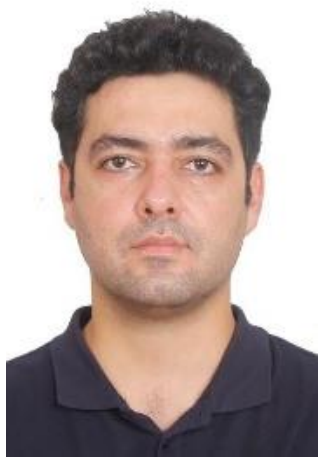

**Hassan rezai Soleymanpour** was born in Karaj, Iran (1985). He received his bachelor's degree in electrical engineering from Hormozgan University, Hormozgan, Iran (2007) and his master’s degree in electrical engineering from Semnan University, Semnan, Iran (2010). His research interests include electrical power system, artificial intelligence, and mathematical problems.